\documentclass[12pt,a4paper,twoside]{amsart}

\usepackage[T1]{fontenc}
\input xypic\usepackage{amsmath}
\usepackage{amssymb}

\usepackage{amssymb}
\usepackage{latexsym}
\usepackage{psfig}
\usepackage{graphicx}
\usepackage{psfrag}
\newcommand{\newtext}{} 

\newcommand{\R}{{\mathbb R}}

\newcommand{\Z}{{\mathbb Z}}

\newcommand{\SO}{{Schr\"odinger operator }}
\newcommand{\ip}{{inverse problem }}
\newcommand{\IBSP}{{inverse boundary spectral problem }}
\newcommand{\BSD}{{boundary spectral data }}
\newcommand{\Rm}{{Riemannian manifold }}

\newcommand{\csd}{{local spectral data }}
\newcommand{\SOs}{{Schr\"odinger operators }}
\newcommand{\wrt}{{with respect to }}
\newcommand{\wlg}{{without loss of generality, }}

\def\hat{\widehat}
\def\tilde{\widetilde}
\def \bfo {\begin {eqnarray*} }
\def \efo {\end {eqnarray*} }
\def \ba {\begin {eqnarray*} }
\def \ea {\end {eqnarray*} }
\def \beq {\begin {eqnarray}}
\def \eeq {\end {eqnarray}}

\def \e {\varepsilon}
\def \p {\partial}

\newcommand{\la}{{\lambda}}

\newtheorem{definition}{Definition}[section]
\newtheorem{theorem}[definition]{Theorem}
\newtheorem{lemma}[definition]{Lemma}
\newtheorem{problem}[definition]{Problem}

\newtheorem{corollary}[definition]{Corollary}
\newtheorem{remark}[definition]{Remark}

\begin{document}
\title{Multidimensional Borg-Levinson theorem}
\author{Yaroslav Kurylev}
\address{Yaroslav Kurylev, Department of Mathematical Sciences, Loughborough Univ.,
Loughborough, LE11 3TU, UK}
\email{Y.V.Kurylev@lboro.ac.uk}

\author{Matti Lassas}\address{Matti Lassas, Helsinki University of
Technology, Institute of Mathematics, P.O.Box 1100,  02015, Finland }
\email{mjlassas@math.hut.fi}
\author{Ricardo Weder}\address{Ricardo Weder, Instituto de Investigaciones en Matem\'aticas
Aplicadas y en Sistemas, Universidad  Nacional Aut\'onoma de
M\'exico, Apartado Postal 20-726, M\'exico DF 01000, M\'exico.}
\email{weder@servidor.unam.mx.}

\maketitle {\bf Abstract.} {\it We consider the inverse problem of
the reconstruction of a Schr\"odinger operator on a unknown
Riemannian manifold or a domain of Euclidean space. The data used
is a part of the boundary $\Gamma$ and the eigenvalues
corresponding to a set of impedances in the Robin boundary
condition which vary on $\Gamma$. The proof is based on the
analysis of the behaviour of the eigenfunctions on the boundary as
well as in perturbation theory of eigenvalues. This reduces the
problem to an inverse boundary spectral problem solved by the
boundary control method. }


\section{Introduction}

 In 1929  Ambartsumyan \cite{Am} considered the
Sturm-Liouville problem
\beq
\label{1}-\psi''+q(x)\,\psi=\lambda\,\psi, \quad x\in (0,1),
\quad
\psi'(0)=\psi'(1)=0,
\eeq
where
the
potential $q$ is continuous and real valued. Let
$\{\lambda_k\}_{k=0}^\infty$ be the eigenvalues  for this
Sturm-Liouville problem. Ambartsumyan proved that if
$\lambda_k=k^2$ for $k=0,1,\dots,$ then $q\equiv 0.$

The next important contribution was due to {Borg} \cite{Bo}
who assumed that $q$ is integrable and real valued. His result can
be stated as follows. He proved that one spectrum in general does
not uniquely determine the corresponding Sturm-Liouville operator
and that the result of Ambartsumyan is a special case.

{\newtext Let} $\{\lambda_k\}_{k=0}^\infty$ be the eigenvalues for (\ref{1})
with the boundary conditions \bfo
 \psi'(0)+ h_1\cdot\psi(0)=0,\quad
\psi'(1)+ h_3\cdot\psi(1)=0,
\efo
and let
$\{\mu_k\}_{k=0}^{\infty}$ be the eigenvalues with the boundary
condition \bfo \psi'(0)+h_2\cdot\psi(0)=0,\quad \psi'(1)
+h_3\cdot\psi(1)=0, \efo where $h_1 \neq h_2, h_3$ are real
numbers. Then, the two sets $\{\lambda_k\}_{k=0}^\infty$ and
$\{\mu_k\}_{k=0}^{\infty}$ uniquely determine $h_1, h_2, h_3$ and
$q$. Levinson \cite{Lev} obtained simpler proofs of some of the
results of Borg.

Borg \cite{Bo1} and Marchenko \cite{Ma} generalized the
Borg-Levinson theorem to Sturm-Liouville operators on the half
line with a boundary condition at the origin when there is no
continuous spectrum. They independently proved  that the discrete
spectra corresponding to two different boundary conditions at
$x=0$ (with a fixed boundary condition, if required, at
$x=+\infty$) uniquely determine the potential and the boundary
conditions at the origin.

Borg-Marchenko's result was generalized to the case where
there is also a continuous spectrum in \cite{AW}  where it was
proven that the potential and boundary conditions are uniquely
determined by an appropriate data set containing the discrete
eigenvalues and continuous part of the spectral measure
corresponding to one boundary condition at the origin and a subset
of the discrete eigenvalues for a different boundary condition.
Another extension of the Borg-Marchenko theorem to the case
with a continuous spectrum is given by Gesztesy and Simon
\cite{GS}. The uniqueness result is proven there  in  the case
when Krein's spectral shift function is known.

The Borg-Levinson inverse two spectra problem can be reduced
to the inverse boundary spectral problem with data of the form
\beq \label{2} \{ \lambda_k, c_k \}_{k=0}^\infty \ \eeq where
$c_k$ are the norming constants, \bfo
  c_k:= \|\psi_k\|_{L^2(0,1)},
\efo
and $\psi_k$ is the eigenfunction corresponding to $\lambda_k$ with
$\psi_k(0)=1, \psi'_k(0)=-h_1$. See for example \cite{Le}, \cite{GL}.
Clearly, data (\ref{2}) is equivalent to the following inverse boundary
spectral data,
\beq
\label{9.5.1}
\{
 \lambda_k, \phi_k(0)\}_{k=0}^\infty,
\eeq
 where now $\phi_k$ are the unit-norm eigenfunctions.

 A multidimensional analog of boundary spectral data is the set
\bfo \{\lambda_k, \phi_k|_{\partial\Omega}\}_{k=0}^\infty, \efo in
the case of the Neumann or third-type boundary conditions (cf.
(\ref{9.5.1}), and the  set \bfo \{\lambda_k,
\partial_n\phi_k|_{\partial\Omega}\}_{k=0}^\infty, \efo in the
case of the Dirichlet boundary {condition}. Here $\Omega
\subset \R^n$ is a (smooth) bounded domain and $\p_n$ is the
interior unit
normal derivative to $\p \Omega$. In comparison with the
$1-$dimensional case, not all second-order elliptic operators,
even isotropic ones, can be reduced to a {Schr\"odinger}
operator in $\Omega$. For different classes of isotropic elliptic
operators, e.g. for
 an
acoustic operator, or a \SO, or a more general second-order
operator, namely, \beq \label{5} \\
\nonumber Au= -c^{-2}(x) \Delta u, \,\,
\hbox{or} \,\, Au= -\Delta u+q(x)u, \,\, \hbox{or} \,\, Au= -
\hbox{div} \left(\e(x) \triangledown u \right)+q(x)u, \eeq where
$c,\, \e$ are positive functions and $q$ is a real-valued function
in $\Omega$, the uniqueness of determination of $c$, or $q$, or
$\e$ and $q$ was proven, correspondingly in \cite{Be}, \cite{NSU}
and \cite{Nov}.
It should be noted that the methods used in these
papers differed {significantly,} with \cite{Be} introducing
the boundary control (BC) method while \cite{NSU} being based on
the complex geometric optics method of \cite{SU} and \cite{Nov}
using the ideas of $\overline{\p}$-problem.

The \IBSP for the
anisotropic case was considered in \cite{BK}, where it was shown
that \BSD determine a compact \Rm and in \cite{Ku1}, \cite{Ku2}
and \cite{KL} where it was shown that \BSD determine,  up to a
natural group of  gauge transformations,  a general second-order
self-adjoint elliptic operator and a
 wide class of second-order non-self-adjoint elliptic operators on a compact
manifold.
It should be noted that, the boundary $\p \Omega$ of the manifold
being given, the manifold itself was not a priori known and was to be
recovered from the \BSD which, in this case, is the set
\beq
\label{BSD1}
\left(\p \Omega,\,\, \{\la_k, \, \phi_k|_{\p \Omega} \}_{k=1}^{\infty} \right)
\eeq
where $\la_k$ and $\phi_k$ are the Neumann-eigenvalues
and normalized eigenfunctions of the Laplace-Beltrami operator.

{\newtext In this paper we use
invariant formulation of inverse problems, i.e., formulate the
problem in terms of manifolds. For clarity, we also apply the obtained
results in the Euclidean setting.
Unless otherwise specified,
$(\Omega,g)$ is a smooth connected compact Riemannian manifold with
non-empty boundary. On $(\Omega,g)$ we study
the Schr\"odinger operator
\ba
A=-\Delta+q
\ea
where $\Delta=\Delta_g$ is the Laplace-Beltrami operator.
By $A^\omega$ we denote the operator $A$ defined in the set of
$H^2(M)$ functions that satisfy the third-type boundary condition on $\p \Omega$,
\bfo
 \left(\p_{\nu} u+\omega u\right)|_{\p \Omega} =0,
\efo with $\p_{\nu}$ being the interior normal derivative on $\p \Omega$ in
the corresponding metric. Following physical literature, we
refer to the real valued function $\omega\in C^\infty(\p \Omega)$ as the impedance.}
The proofs  in \cite{Ku1},
\cite{Ku2}, \cite{KL} were based on a geometric approach to the
BC-method, see \cite{KKL} for a detailed exposition. It is,
however, clear from the considerations above that the mentioned papers
on multidimensional inverse problems did not consider a multidimensional analog
of the Borg-Levinson \ip, but  the inverse boundary spectral problem.
A
multidimensional analog of the Borg-Levinson inverse problem
may be formulated as follows:

\begin{definition}
\label{CSD}
Let $(\Omega,g)$ be a compact connected \Rm with non-empty
boundary $\p \Omega$, $\Sigma \subset \p \Omega$ be an open connected
non-empty subset  and $q $ be a real-valued function in $ C^{\infty}(\Omega)$.
Let $\omega_0 \in C^{\infty}(\p \Omega)$ be a real valued function. Consider
 the
\SOs in $L^2(\Omega)$ of the form, \beq \label{SO} A^{\omega}u =
-\Delta u +qu, \quad D(A^{\omega}) = \{u \in H^2(\Omega):
\left(\p_{\nu} u+\omega u\right)|_{\p \Omega}  =0\},
\eeq
where $\omega$ is real valued and
$\tilde{\omega}= \omega-\omega_0 \in C^{\infty}_0(\Sigma)$.
{\newtext Denote
by $\la_k(\omega),\,k=1, 2, \dots$ the corresponding eigenvalues
counting multiplicity. The \csd is }
 \beq \label{BLSD} \hbox{ }\ \ \
\Sigma\ \hbox{and the map $\omega\mapsto
\{\la_k(\omega) \}_{k=1}^{\infty}$ defined for
$\omega \in C^{\infty}(\p \Omega)$, $\omega-\omega_0 \in
C^{\infty}_0(\Sigma)$}. \eeq
\end{definition}

Note that here $\Omega$ is compact manifold so that
 $ C^{\infty}(\Omega)$ consists of functions that are smooth upto
the boundary.

\begin{problem}
\label{BLSP}
Do \csd of form (\ref{BLSD}) determine $(\Omega,g)$, $q$ and $\omega_0$ uniquely?
\end{problem}

\noindent Note, that by  determination of a Riemannian manifold $(\Omega,g)$ we
mean determination
of its isometry type.

{\newtext
We denote the Gateaux derivatives of
$\omega\mapsto \la_k(\omega)$ at $\omega_0$ in the direction $\tilde{\omega}$
by $\la_{k, \omega_0}(\tilde{\omega})=d\la_k|_{\omega_0}(\tilde \omega).$
Clearly, \csd make it possible to find the
 $ \la_{k, \omega_0}(\tilde{\omega})$ for any
$k=1,2,\dots$ and $\tilde \omega\in C^{\infty}_0(\Sigma)$.

 In following, we
use notation
\beq \label{ball}
B_{\e}^{\infty}(\omega_0) = \{\omega\in C^\infty(\p\Omega): \ ||\omega-
\omega_0||_{L^{\infty}(\p \Omega)} < \e,\
 \omega - \omega_0
\in C^{\infty}_0(\Sigma)\}.
\eeq
}
Depending on degeneracy/non-degeneracy of the
spectrum of $A^{\omega_0}$, we prove the following result.
\begin{theorem}
\label{local} Let {\newtext $(\Omega,g)$ be a
smooth, compact, connected Riemannian manifold with boundary and}
$\Sigma \subset \p \Omega$ be an open,
connected, non-empty subset and $A^{\omega_0}$ be
a \SO of form (\ref{SO}). Then

\smallskip \noindent
{\it a.}  If the spectrum of $A^{\omega_0}$ is simple, then
$\Sigma$, the eigenvalues $\la_k(\omega_0)$, and their
Gateaux derivatives, $\la_{k, \omega_0}(\omega)$; $
\omega \in C_0^{\infty}(\Sigma)$ uniquely determine $(\Omega,g),\,
q$ and $ \omega_0$.

\smallskip \noindent
{\it b.} {\newtext  For arbitrary $A^{\omega_0}$,} given $\Sigma$ and
$\{\la_k(\omega)\}_{k=1}^{\infty}$ for all real-valued $\omega \in
B^{\infty}_{\epsilon}(\omega_0)$
with  some $\epsilon >0$, one can uniquely determine
 $(\Omega,g),\, q$ and $ \omega_0$.
\end{theorem}

\noindent Note that, in Theorem \ref{local}, we do not assume an a
priori knowledge of either $\Omega$ or $\p \Omega$.  We only
have to know $\Sigma$.
{\newtext Theorem \ref{local} has the following corollary in Euclidean
setting.}

\begin{corollary}\label{Euclidean} {\newtext Let $\Omega\subset \R^n$,
$g_{ij}(x)=c(x)\delta_{ij}$ be a conformally isotropic metric on $\Omega$,
and $\Sigma\subset \p \Omega$ be open and non-empty.
Let $A^{\omega_0}$
be a \SO of form (\ref{SO}).
Then $\Sigma$ and
$\{\la_k(\omega)\}_{k=1}^{\infty}$ for all real valued $\omega \in
B^{\infty}_{\epsilon}(\omega_0)$,
with  some $\epsilon >0$, determine $\Omega$
as a subset of $\R^n$, $c(x)$, $q$, and $ \omega_0$ uniquely.
}
\end{corollary}

\section{Boundary behavior of eigenfunctions}
In this section we consider the eigenvalues and eigenfunctions of an operator $ A^{\omega}$
for a fixed $\omega$. In this connection we skip using $\omega$ throughout this
section, writing  $\la_k$ instead of $\la_k(\omega)$ and
$\phi_k$ instead of $\phi_k(\omega)$.

To describe behavior of {\newtext eigenfunctions} near $\p \Omega$ we employ the boundary normal
coordinates $x=(z,\tau)$, where $\tau =\hbox{dist}(x,\p \Omega)$ and $ z$ is the unique
point  on $\p \Omega$ nearest to $x$ with local coordinates $z=(z^1,\dots,z^{n-1})$.

\begin{lemma}
\label{unique cont}
Let $\phi$ be an eigenfunction for an eigenvalue $\la$ of an operator $ A^{\omega}$
(with some fixed $\omega$). Then, for any $z_0 \in \p \Omega $, there is a
multi-index
$\alpha_0 \in \Z_+^{n-1}$ such that
\beq
\label{11.1}
\p^{\alpha_0} \phi(z_0) \neq 0.
\eeq
\end{lemma}
Here $\phi(z)=\phi(z,0)$ and equation (\ref{11.1}) is valid in proper local coordinates
on $\p \Omega$, $z=(z^1,\dots,z^{n-1})$ where, \wlg $z_0=0$.

{\bf Proof.} If $\omega \neq 0$ we introduce a gauge
transformation {\bf \cite{KKL}} \bfo u \longrightarrow v = \kappa
u, \quad \kappa \in C^{\infty}(\overline{\Omega}), \quad \kappa(x)
>0 \,\, \hbox{for} \,\, x \in \overline{\Omega}, \,\, \p_{\tau}
\kappa|_{\tau =0} = - \omega. \efo Then $\psi = \kappa \phi$ is a
smooth solution to the equation \beq \label{11.2}
\\ \nonumber -
 \p^2_{\tau} \psi - g^{ij} \p_i \p_j \psi + a^n
\p_{\tau} \psi + a^i \p_i \psi +a^0 \psi =\la \psi, \,\, \tau
>0,\,\, i,j=1,\dots,n-1, \eeq where $a^0, a^i,  a^n$ and $g^{ij}$
are functions of $(z,\tau)$, and \beq \label{11.3} \p_{\tau}
\psi|_{\tau=0}=0. \eeq
{\newtext Assume now that, for any $\alpha=
(\alpha_1,\dots,\alpha_{n-1}) \in
\Z_+^{n-1}$, $
\p^\alpha  \phi(0)=
\p_{z_1}^{\alpha_1}\dots\p_{z_{n-1}}^{\alpha_{n-1}}  \phi(0)=0$ and, therefore, $\p^{\alpha}
\psi(0) =0$. Using (\ref{11.2}), (\ref{11.3}),  this implies that
for any $\beta=(\beta_1,\dots,\beta_n) \in \Z_+^n$,
\beq \label{11.4}
\p_{z_1}^{\beta_1}\dots\p_{z_{n-1}}^{\beta_{n-1}}\p_\tau
^{\beta_{n}} \psi(0)
=0. \eeq
}

Let $\hat{\psi}, \hat{a}^0, \hat{a}^i, \hat{g}^{ij}$ be even
continuations of these functions across the boundary $\tau=0$ and
$\hat{a}^n$ be an odd continuation of $a_n$. Then, in an open set
$U \subset \R^n, 0 \in U$, the function $\hat{\psi}$ is a $C^2(U)$
solution of the equation \beq \label{11.5} -\p^2_{\tau} \hat{\psi}
- \hat{g}^{ij} \p_i \p_j \hat{\psi} + \hat{a}^n  \p_{\tau}
\hat{\psi} + \hat{a}^i \p_i \hat{\psi} +\hat{a}^0 \hat{\psi} = \la
\hat{\psi}, \eeq with $\hat{g}^{ij}  \in C^{0,1}(U)$ and
$\hat{a}^p \in L^{\infty}(U),\, p=0,\dots,n$. Moreover, by
(\ref{11.4}), for any $N>0$ there is $C_N$ so that \bfo
|\hat{\psi}(z,\tau)| \leq C_N |x|^N, \quad |x|^2 =
\sum_{i=1}^{n-1} |z^i|^2+ \tau^2. \efo This, together with
equation (\ref{11.5}) imply, due to the H\"ormander strong
uniqueness principle, \cite{Ho}, that $\psi=\phi =0$. \hfill$\Box$

It will be shown in the next section that, under some additional assumptions,
\csd determine $|\phi_k^{\omega_0}(x)|, \, x \in \Sigma,\, k=1,  2, \dots$. Moreover, the
following result holds:

\begin{theorem}
\label{th:2}
Given $\xi \in C^{\infty}(\Sigma)$ such that  $\xi(z) = |\phi(z)| , \, z \in \Sigma$, where $\phi$ is an eigenfunction of an
operator $A^{\omega_0}$, it is possible to find  $ \phi|_{\Sigma}$ up to multiplication by $\pm 1$.
\end{theorem}

{\bf Proof.} To fix the sign of $\phi$, choose a point  $z_0
\in \Sigma$ where $\xi(z_0) >0$ and take
$\phi(z_0) =\xi(z_0) >0$. Let $(z^1,\dots,z^{n-1}) \in B_r
\subset \R^{n-1}$ be Riemannian normal coordinates
in the metric ball $B_r(z_0) \subset \Sigma$,
{\newtext where $(\p \Omega,g)$ is endowed with the metric induced by
$(\Omega,g)$}.
Note that we can choose
\beq
\label{27.1}
r= \min (\hbox{inj}(\p \Omega), \, d_{\p \Omega}(z_0, \p \Sigma)).
\eeq

We first show that $\xi$ determines $\phi$ everywhere in $B_r$. By
continuity of $\phi$, it is clear that $\xi$ determines $\phi$ in
ball $B_{\tilde{r}}$ for sufficiently small $\tilde{r}$.
 {\newtext Let $\rho$ be the largest possible
value $\rho\leq r$ such that
 $\xi$ determines $\phi$ in
ball $B_{\rho}$.}
We want to show that $\rho=r$. Assuming the
contrary, we note that $\phi$ is defined on the closure
$\overline{B_{\rho}}$. Let $y \in \p B_{\rho} \subset \Sigma$. If
$\phi(y) \neq 0$,  then, by continuity, $\xi$ determines $\phi$ in
a vicinity of $y$. If $\phi(y) = 0$, by Lemma \ref{unique cont}
there is $m>0$ such that
\beq
\label{11.7} \phi(z)= \sum_{|\alpha|
=m} b_{\alpha}(z-y)^{\alpha} +O\left(|z-y|^{m+1}\right), \quad
b_{\alpha_0} \neq 0 \,\, \hbox {for some} \,\, \alpha_0,\
|\alpha_0|=m.
\eeq
It follows from (\ref{11.7}) that there is an open dense set $W
\subset S^{n-2}$ such that, for $e = (e^1,\dots,e^{n-1}) \in W$,
\bfo \p^m_e \phi (y) = (e^j \p_j)^m \phi(y) \neq 0. \efo Choosing
$e$ transversal to $\p B_{\rho}$ at $y$ and assuming that, \wlg
$\p_e = \p_{n-1}$, we obtain, using Malgrange Preparation Theorem,
e.g. \cite[ Th. 7.5.5]{Ho1}, that, in a vicinity of $y$, \beq
\label{malgrange} \phi(\hat{z}, z^{n-1}) = c(\hat{z}, z^{n-1})
\, \sum_{l=0}^m a_l(\hat{z}) (z^{n-1}-y^{n-1})^l. \eeq Here
$\hat{z} =(z^1,\dots, z^{n-2})$, the function $c(z)= c(\hat{z},
z^{n-1})$ is a $C^{\infty}-$function near $z=y$ with $ c(y)\neq
0$ and $a_m(\hat{z}) =1$. Therefore, $\phi(\hat{z}, \,z^{n-1})$,
for a fixed $\hat{z},$ has only a finite number of real roots,
$r_j(\hat{z})$. The function $\phi(\hat{z}, z^{n-1})$, considered
as a function of $z^{n-1}$, changes its sign at $r_j(\hat{z})$
when this root is of an odd order and does not change the sign
when the root is of an even order. As the lines $\hat{z} =
\hbox{const}$ are transversal to $\p B_{\rho}$ near $y$, we obtain
the continuation of $\phi$ into a vicinity of $y$. As $y \in \p
B_{\rho}$ is arbitrary, we obtain the continuation of $\phi$ into
an open neighborhood of $\overline{B_{\rho}}$.  Thus, $\rho =r$,
i.e. $\phi$ can be uniquely determined everywhere in the ball
$B_{r}$. It also follows from the above arguments that $\{z:
\phi(z) \neq 0\} \cap B_r(z_0)$ is an open set of full measure.

To proceed further, let $\tilde{z} \in \Sigma$ and
 $L$ {\newtext be} a curve in $\Sigma$ connecting $z_0$ with $\tilde{z}$.
We cover $L$ by a finite number of balls $B_{r_j/2}(z_j), \, j=0,1,\dots, J,\,
z_J=\tilde{z},$ such that $ B_{r_j/2}(z_j) \cap B_{r_{j+1}/2}(z_{j+1}) \neq \emptyset$.
In particular, there is a point $\tilde{z_1} \in B_{r_0/2}(z_0) \cap B_{r_1/2}(z_{1})$
with $\phi(\tilde{z}_1) \neq 0$.  By the previous construction we find $\phi$ in
$B_{\tilde{r}_1}(\tilde{z}_1)$  which contains $ B_{r_1/2}(z_{1})$. Continuing this process, we
find $\phi(\tilde{z})$.
\hfill$\Box$

\section{Generic behavior of eigenvalues}

Consider the quadratic form $Q^{\omega}$ related to the operator $A^{\omega}$,
\beq
\label{quadratic form}
Q^{\omega} (u) = \int_{\Omega} \left(| \nabla u|^2 +q|u|^2 \right)\,dV
+\int_{\p \Omega} \omega |u|^2 \,dS,
\eeq
where $dV,\, dS$ are the volume and area forms generated by the metric $g$ in $\Omega$
and $\p \Omega$.

Let $A(t)$ be an  analytic, for $|t| < \e$, one-parameter family of \SOs of the form (\ref{SO}), where
the impedance $\omega(t)$ of the form
\beq
\label{1-family}
\omega(t)=\omega_0+t \tilde{\omega}, \quad \hbox{with real}\,\,\tilde{\omega}
\in C^{\infty}_0(\Sigma).
\eeq
 Then $A(t)$ is a self-adjoint homomorphic operator family
of type (B), in the sense of Kato \cite[ Section 7.4]{Ka},
so that the eigenvalues $\la_k(\omega(t))$ and eigenfunctions
$\phi_k^{\omega(t)}$ may be chosen to be analytic with respect to $t$.
 In this case we can find the Gateaux derivative of
$\la_k$ \wrt $t$. A bit more generally, the following result holds:

\begin{lemma}
\label{lem.1}
 Let  $\la_k(t),\, \phi_k(t)$ be an eigenvalue and a corresponding normalized
eigenfunction of $A(t)$ which are differentiable \wrt  $t$.
Then
\beq
\label{gateau}
\dot{\la}_k(t)= - \int_{\p \Omega} |\phi_k(t)|^2 \tilde{\omega} \, dS,
\eeq
where $\dot{\la}$ stands for the $t-$differentiation of $\la$.
\end{lemma}

{\bf Proof.}
Differentiating \wrt $t$ the equation for $\phi_k(t)$, we get
\bfo
\left(-\Delta +q -\la_k(t) \right) \dot{\phi}_k(t)= \dot{\la}_k(t) \phi_k(t).
\efo
Thus, due to $||\phi_k(t)||=1$,
\beq
\label{10.1}
\dot{\la}_k(t)=\int_{\Omega} \left((-\Delta +q -\la_k(t) ) \dot{\phi}_k(t)  \right)\,
\overline{\phi_k(t)} \, dV \\
\nonumber
=\int_{\p \Omega} \left( \p_{\nu} \dot{\phi}_k(t)\, \overline{\phi_k(t)} -
 \dot{\phi}_k(t)\, \p_{\nu} \overline{\phi_k(t)}\right) \,dS.
\eeq
By the boundary condition in (\ref{SO}),
\bfo
\p_{\nu} \dot{\phi}_k(z,t) = - \left(\omega(z, t)\dot{\phi}_k(z,t) +
\dot{\omega}(z,t)\phi_k(z,t)\right), \quad z \in \p \Omega.
\efo
This together with (\ref{10.1}) imply
equation (\ref{gateau})  due to (\ref{1-family}).
\hfill$\Box$

Denote by   $\mu_k(\omega)$  the multiplicity of
$\la_k^{\omega}$ and assume that $\mu_k(\omega)$ is constant near
$\omega_0$.


\begin{corollary}
\label{cor:1} Assume that for some $\e >0$,  $\la_{k-j-1}(\omega)
<\la_{k-j}(\omega)= \dots \la_k(\omega)= \dots
=\la_{k+p-1}(\omega) < \la_{k+p}(\omega)$, $p+j= \mu_k(\omega_0)$,  for
{\newtext all}\beq
\label{11.6} \omega \in B_{\e}^{\infty}(\omega_0). \eeq Then for
any $\tilde \omega\in C^\infty_0(\Sigma)$
and a normalized eigenfunction $\phi$ of $A^{\omega_0}$ corresponding to
the eigenvalue $\la_k(\omega_0)$ there is an eigenvalue $\lambda(t)$
and a normalized eigenfunction
$\phi(t)$ of $A^{\omega(t)}$, $\omega(t)=\omega_0+t\tilde \omega$
 such that $\phi(0)=\phi$ and that the
equation
(\ref{gateau}) is valid.
\end{corollary}

This result is standard for the perturbation theory for quadratic forms, e.g. \cite{Ka}, \cite{BS}.
We repeat its proof for the convenience of the reader.

{\bf Proof.} By the perturbation theory for quadratic forms, e.g. \cite{Ka}, \cite{BS},
a sufficiently small disk centered in $\la_k=\la_k(\omega_0)$  does not contain eigenvalues of
$A^{\omega}$, except for $\la_{k-j}(\omega), \dots, \la_{k+p-1}(\omega)$,  when $\omega$ satisfies (\ref{11.6}) with
sufficiently small $\e$.
Consider the Riesz projectors, $P_k^{\omega}$,  to the eigenspace
corresponding to $\la_k(\omega)$,
\beq
\label{riesz}
P_k^{\omega} = \frac{1}{2 \pi i} \int_{\Gamma} R_z^{\omega}\, dz,
\eeq
where $ R_z^{\omega}$ is the resolvent for $A^{\omega}$ and $\Gamma$ is a
sufficiently small circle around $\la_k(\omega_0)$. When $\omega=\omega(t)$ is
of form (\ref{1-family}), $ R_z(t)$  is an analytic, \wrt $t$,
operator-valued function in $L^2(\Omega)$. Therefore,
$ P_k^{\omega}$
are also analytic \wrt $t$. Moreover, for sufficiently small $\e$ and real $t$,
$\tilde{\phi}(t) = P_k^{\omega(t)} \phi  \neq 0$ so that
$\phi(t) =  \tilde{\phi}(t) /||\tilde{\phi}(t) ||$ is a desired  normalized
eigenfunction for $A(t)$ which smoothly depends on $t$.  This implies also that
$\la_k(t)$ is smooth \wrt $t$ and the considerations of
Lemma \ref{lem.1} are valid.
\hfill$\Box$

Combining Corollary \ref{cor:1} with Theorem \ref{th:2} we obtain the following
result.

\begin{corollary}
\label{cor:2}
Assume that $\la_k(\omega)$ has a constant multiplicity, $\mu_k(\omega)=\mu_k(\omega_0)$
for all $\omega$ satisfying
equation (\ref{11.6}). Then $\mu_k(\omega)=1$.
\end{corollary}

{\bf Proof.}
By corollary \ref{cor:1}, any $\phi \in P_k^{\omega_0} L^2(\Omega), \,
||\phi||=1$ satisfies equation (\ref{gateau}). Thus, for any two different normalized
eigenfunctions $\phi,\, \tilde{\phi}$ for $\la_k$,
\bfo
\int_{\p \Omega}  |\phi|^2 \tilde{\omega}\,dS =
\int_{\p \Omega}  |\tilde{\phi}|^2\tilde{\omega} \,dS,
\efo
with arbitrary $\tilde{\omega} \in C^{\infty}_0(\Sigma)$.This implies that
$|\phi| =|\tilde{\phi}|$ on $\Sigma$, so that
 $\phi|_{\Sigma} = \pm \tilde{\phi}|_{\Sigma}$. This, together
with the boundary condition in (\ref{SO}), yield that also
$\p_{\nu} \phi|_{\Sigma} = \pm \p_{\nu}\tilde{\phi}|_{\Sigma}$.
Using the similar arguments as in proof of Lemma \ref{unique cont} and applying the
H\"ormander unique continuation theorem \cite{Ho}, $\phi= \pm \tilde{\phi}$ on $\Omega$.
\hfill$\Box$

We now investigate the multiplicity of eigenvalues under small perturbations of the
impedance.

\begin{lemma}
\label{lm:3} For any $k \in \Z_+,\, \e >0$ there is $\omega \in
C^{\infty}(\partial \Omega)$ satisfying equation
(\ref{11.6}) such that $\la_i(\omega)$ are simple for $i=1, \dots,
k$.
\end{lemma}

{\bf Proof.} By the perturbation theory for quadratic forms e.g.
\cite{Ka}, \cite{BS}, for any $\omega_0 \in C^{\infty}(\p \Omega)$
and $i \in \Z_+$, there are $\epsilon, \, \delta >0$ such that
\bfo \hbox{dim} \, P^{\omega}(\delta) L^2(\Omega) = \hbox{dim} \,
P^{\omega_0}_i L^2(\Omega), \efo for all $\omega \in
B^{\infty}_{\epsilon}(\omega_0)$, where $ P^{\omega}(\delta)$ is
the projector onto the sum of eigenspaces of $A^{\omega}$
corresponding to the eigenvalues from the interval $(
\la_i(\omega_0) - \delta, \, \la_i(\omega_0) + \delta)$.
Therefore, $\mu_i(\omega_0)$ is an
upper-semicontinuous function of $\omega_0 \in L^{\infty}(\p
\Omega)$.

Let \bfo \underline{\mu}_i(\omega_0)= \liminf_{\omega \to
\omega_0,\ \omega-\omega_0\in C^\infty_0(\Sigma)}
\mu_i(\omega),\,\, \overline{\mu}_i(\omega_0)=
\limsup_{\omega \to \omega_0,\ \omega-\omega_0\in C^\infty_0(\Sigma)}
 \mu_i(\omega) .
\efo
 As $\mu_i(\omega) \in \Z_+$, there is $\delta_i=\delta_i(\omega_0) >0$, such that
\bfo \min \mu_i(\omega) =\underline{\mu}_i(\omega_0), \,\, \max
\mu_i(\omega) =\overline{\mu}_i(\omega_0), \efo where {\newtext minimum and
maximum} are  taken over the set
$\omega\in B^{\infty}_{\delta_i(\omega_0)}(\omega_0).$
 Choose $\omega_1$ with $\mu_1(\omega_1)=\underline{\mu}_1(\omega_0)$ such that
\bfo
||\omega_1 - \omega_0||_{ L^{\infty}(\p \Omega)} < \min (\e/k,\,
\delta_1(\omega_0)),
\quad \omega_1-\omega_0 \in C^{\infty}_0(\Sigma).
\efo
Then, due to the mentioned upper-semicontinuity of $\mu_1$, there is $\tilde{\delta}_1 >0$
so that
\beq
\label{11.8}
\mu_1(\omega) = \mu_1(\omega_1)
\quad \hbox{for} \,\  \omega\in B^{\infty}_{\tilde \delta_1}
(\omega_1).
\eeq
By Corollary \ref{cor:2},
\bfo
\underline{\mu}_1(\omega_0) = \mu_1(\omega)=1,
\efo
for $\omega\in B^{\infty}_{\tilde \delta_1}(\omega_1).$

Next we find $\delta_2<\tilde{\delta}_1$ such that \bfo \min\mu_2(\omega) =\underline{\mu}_2(\omega_1) \efo where
minimum is taken over the set \beq \label{11.9} \omega \in
B_{\delta_2}^{\infty}(\omega_1).
\eeq
This makes it possible to choose
 $\omega_2$ satisfying (\ref{11.9}) and also
\bfo \mu_2(\omega_2) =\underline{\mu}_2(\omega_1), \quad
||\omega_2-\omega_1||_{L^{\infty}(\p \Omega)} < \min(\e/k,\,
\delta_2). \efo Repeating the same arguments as for $\mu_1$, there
is $\tilde{\delta}_2 < \min(\e/k,\, \delta_2)$ such that \beq
\label{11.10} \underline{\mu}_2(\omega_1) = \mu_2(\omega)=1, \eeq
for $\omega \in  B^{\infty}_{\tilde{\delta}_2}(\omega_2)$,
 and the ball
$B^{\infty}_{\tilde{\delta}_2}(\omega_2)$ lies inside the ball
$B_{\tilde{\delta}_1}^{\infty}(\omega_1)$  so that also \bfo
\mu_1(\omega)=1. \efo Continuing this procedure, we find $\omega_k
\in C^{\infty}(\p \Omega), \omega_k-\omega_0 \in C^{\infty}_0(
\Sigma)$ with \beq \label{11.12} \mu_1(\omega_k)
= \dots =\mu_k(\omega_k) =1. \eeq Moreover, it is seen easily from
the above construction that \bfo ||\omega_k - \omega_0||_{
L^{\infty}(\p \Omega)} < \e. \efo \hfill$\Box$

\begin{remark}
\label{rem:1}
A slight  modification of the previous arguments shows that, in any
$ C^{\infty}_0(\Sigma)-$neighborhood of $\omega_0$ there is
an impedance $\omega$ such that the spectrum of $A^{\omega}$
is simple. Indeed, we can easily generalize Lemma \ref{lm:3} to show that,
for any $k \in \Z_+,\, \e >0$ and $\omega$ there is $\omega_k$ satisfying (\ref{11.12})
such that
\beq
\label{11.13}
|| \omega_k - \omega||_{ C^k(\p \Omega)} < \frac{\e}{2^k}.
\eeq
To construct $\omega$ with simple spectrum, we first find $\omega_1$ with
$\mu_1(\omega_1)=1$ satisfying (\ref{11.13}) with $k=1$ and $\omega_0$
instead of $\omega$. Then we find $\omega_2$ with
$\mu_1(\omega_2)=\mu_2(\omega_2)=1$ and (\ref{11.13})
with $k=2$ and $\omega_1$ instead of $\omega$. By taking, if necessary,
$\omega_2$ to be $L^{\infty}-$closer
 to $\omega_1$, we obtain that
\beq
\label{11.14}
| \la_1(\omega_2) - \la_2(\omega_2) |>
(1/2 -1/2^{2}) | \la_1(\omega_1) - \la_2(\omega_1) |.
\eeq
Next we find $\omega_3$ with
$\mu_1(\omega_3)=\mu_2(\omega_3)=\mu_3(\omega_3)=1$ and (\ref{11.13})
with $k=3$ and $\omega_2$ instead of $\omega$.  By taking, if necessary,
$\omega_3$ to be $L^{\infty}-$closer
 to $\omega_2$, we obtain that
\beq
\label{11.15}
| \la_1(\omega_3) - \la_2(\omega_3) |>
(1/2-1/2^{3}) | \la_1(\omega_1) - \la_2(\omega_1) |, \\
\nonumber
| \la_2(\omega_3) - \la_3(\omega_3) |>
(1/2-1/2^{3}) | \la_2(\omega_2) - \la_3(\omega_2) |.
\eeq
Continuing the above procedure, we construct a converging, in $ C^{\infty}(\p \Omega)$,
sequence $\omega_k$. Denote by $\omega$ its limit,
$
\omega= \lim \omega_k.
$
By (\ref{11.13}), for any $p \in \Z_+$,
\beq
\label{11.16}
|| \omega_0 - \omega||_{ C^p(\p \Omega)} < \e.
\eeq
As $\la_i(\omega)$ depends continuously on $\omega$, equations
(\ref{11.14}), (\ref{11.15}), and analogous equations for further
$\omega_k$ show that \bfo | \la_k(\omega) - \la_{k+1}(\omega) |
\geq \frac12 | \la_k(\omega_k) - \la_{k+1}(\omega_k) | >0, \efo so
that $A^{\omega}$ has simple spectrum. It is clear from the above
construction that the set of impedances $\omega$ with degenerate
spectrum is of the first Baire category.
\end{remark}
We note that the above result can be also obtained using \cite{Te}, however,
the method of \cite{Te} is different from the one in Remark \ref{rem:1} being based
on the ideas of \cite{Ul} rather than the quadratic forms perturbation theory and unique
continuation for elliptic equation.

\section{From \csd to \BSD. Proof of  main results.}

We are now in the position to prove our main results.  We start with the following technical theorem:

\begin{theorem}
\label{th:3} For any real $\omega_0 \in C^{\infty}(\p \Omega)$ and
any open, non-empty connected  $\Sigma \subset \p \Omega$, the  \csd
determine the traces $\phi_k|_{\Sigma}, \, k=1, \dots,$ {up to
a sign}, of the eigenfunctions of the \SO $A^{\omega_0}$.
\end{theorem}

{\bf Proof.}
If $\mu_i(\omega_0)=1$, Corollary \ref{cor:1} makes possible to find,
for an arbitrary
$\tilde{\omega} \in C^{\infty}_0(\Sigma)$,
\beq
\label{12.3}
\int_{\p \Omega}  \tilde{\omega} |\phi_i|^2 \, dS,
\eeq
where
$\phi_i$ is the normalized eigenfunction of $A^{\omega_0}$ corresponding
to $\la_i(\omega_0)$.

Let now $\mu_i(\omega_0)=p>1$, say $\la_l =\dots=\la_i = \dots=\la_m,\,
l \leq i \leq m,\, m-l=p-1$. By Lemma \ref{lm:3}, there are smooth impedances
$\omega_n, \, n=1,2,\dots,$ which converge to $\omega_0$ while their eigenvalues
$\la_j(\omega_n),\, 1 \leq j \leq m$, remain simple.  By Corollary \ref{cor:1}
it is possible to find $\int_{\p \Omega}  \tilde{\omega} |\phi_j^n|^2 \, dS$,
where
$\phi_j^n,$ for $ l \leq j \leq m$,  are the orthonormalized eigenfunctions
of $A^{\omega_n}$ corresponding
to $\la_j(\omega_n)$.  As $||\phi_j^n||_{H^1(\Omega)}$ are uniformly bounded,
there is a subsequence $n(k)$, which we assume to be the whole
sequence, such that
\beq
\label{12.1}
\lim_{n \to \infty} \phi_j^n = \phi_j, \quad 1 \leq j \leq m.
\eeq
The convergence in (\ref{12.1}) is weak in $H^1(\Omega)$ and strong in
$H^s(\Omega)$ for any $s <1$. As
\bfo
\lim_{n \to \infty} \la_j(\omega_n) = \la_j(\omega_0),  \quad 1 \leq j \leq m,
\efo
$\phi_j$ satisfy the equation
\bfo
(-\Delta+ q) \phi_j = \la_j(\omega_0) \phi_j.
\efo
Moreover,  as
\bfo
\lim_{n \to \infty} \phi_j^n|_{\p \Omega} = \phi_j|_{\p \Omega}
\quad \hbox{in} \,\, L^2(\p \Omega),
\efo
we see that $\phi_j$
are normalized eigenfunctions of $A^{\omega_0}$ for $\la_j, \, 1 \leq i \leq m$.
In addition, for multiple eigenvalues of $A^{\omega_0}$, the corresponding
eigenfunctions remain orthogonal because the eigenfunctions $\phi_j^n,  \phi_k^n,
\,j, k\leq m$ are orthogonal for any $n$ and $j \neq k$.
Thus, $\phi_j$ are the first $m$ orthonormal eigenfunctions of $A^{\omega_0}$.

Also,
\beq
\label{12.2}
\lim_{n \to \infty} \int_{\p \Omega}   |\phi_j^n|^2\tilde{\omega} \, dS =
\int_{\p \Omega}  |\phi_j|^2 \tilde{\omega}\, dS,
\eeq
for any $\tilde{\omega} \in C^{\infty}_0(\Sigma)$, so that
 we know all integrals (\ref{12.3}) when $i \leq m$. Since $m \in \Z_+$ is arbitrary,
we determine the integrals (\ref{12.3}) for any $i \in \Z_+$ and
$\tilde{\omega} \in C^{\infty}_0(\Sigma)$. In turn, this
determines all functions $|\phi_i|_{\Sigma}$. Applying Theorem
\ref {th:2} we find  $\phi_i, \, i=1, 2, \dots,$  on $\Sigma$
up to a sign. \hfill$\Box$


{\bf Proof of Theorem \ref{local}.}

\noindent {\it a.} If all eigenvalues of $A^{\omega_0}$ are
simple, then, by upper semicontinuity of $\mu_k$, it follows from
{Corollary} \ref{cor:1} that the Gateaux derivatives of
$\la_k(\omega_0)$ determine the integrals (\ref{12.3}) for any
$\tilde{\omega} \in C^{\infty}_0(\Sigma)$.  It then follows from
the proof of Theorem \ref{th:3} that the Gateaux derivatives of
$\la_k(\omega_0)$ determine $\phi_i(\omega_0), \, i=1, 2, \dots,$
on $\Sigma$.

\noindent {\it b.}
In  general, Theorem \ref{th:3} shows that $\la_i(\omega)$ for $\omega$ satisfying
(\ref{11.6}) with any $\e >0$, determine $\phi_i(\omega_0), \, i=1, 2, \dots,$
on $\Sigma$.
{\newtext By \cite[Thm. 7.3]{KK}, this data determines uniquely the isometry type
of $(\Omega,g)$ and the gauge-equivalence class
$\{\kappa^{-1}A^{\omega_0}\kappa:\ \kappa\in C^\infty(\Omega),\
\kappa(x)>0\}$ of the operator $A^{\omega_0}$.
By \cite[Lemma 2.29]{KKL} this equivalence class contains a unique Schr\"odinger operator
of the form (\ref{SO}). Thus we can find $q$ and $\omega_0$. This completes
the proof of Theorem \ref{local}. \hfill$\Box$

{\newtext Corollary \ref{Euclidean} is a direct consequence of
 Theorem \ref{local} and the fact that by  Liouville Theorem
\cite{IM}, an isometric
embedding of a conformally Euclidean
$n$-manifold to $\R^n$ is unique.
}

\begin{remark}
 In the case where some of the eigenvalues of $A^{\omega_0}$ are
simple and some are degenerate, our proof gives a result that is
slightly more general than (b) in Theorem 1.3.

We have actually proven that $(\Omega, g), q$ and $\omega_0$ are
uniquely determined by the  data consisting of $\Sigma$,
the simple eigenvalues $\la_k(\omega_0)$ and their
Gateaux derivatives, $  \la_{k, \omega_0}(\omega);$ $ \omega \in
C_0^{\infty}(\Sigma)$, and moreover,
for each degenerate eigenvalue, $\la_k(\omega_0),$ with
multiplicity $\mu_k(\omega_0),$ the local spectral data,
$\{\la_l(\omega)\}_{ l=k-j}^{k+p-1}$ for all $\omega \in
B^{\infty}_{\epsilon}(\omega_0)$
 for some $\epsilon > 0$, and where $p,j$
are the only integers such that, $p+j= \mu_k(\omega_0)$, and
$$
\la_{k-j-1}(\omega_0) < \la_{k-j}(\omega_0)=\cdots=
\la_{k}(\omega_0))= \cdots= \la_{k+p-1}(\omega_0) <
\la_{k+p}(\omega_0).
$$

\end{remark}

{\newtext {\bf Acknowledgements:} Part of the work was done when
the authors were visiting Helsinki within the framework of the
Finnish Inverse Problems Theme Year 2004 R. Weder is a fellow of
Sistema Nacional de Investigadores, was partially supported by
project PAPIIT-UNAM IN 101902, and he thanks M. Lassas and L.
P\"aiv\"arinta for their kind hospitality at Helsinki. Y.~Kurylev
and M.~Lassas were partially supported by the Royal Society. Also,
M.~Lassas's work was supported by Academy of Finland project
102175 and Y.~Kurylev by EPSRC GR/R935821/01. }

\end{document}